\documentclass{article}

\usepackage[english]{babel}
\usepackage[latin1,applemac]{inputenc}
\usepackage{amsmath}
\usepackage{amssymb}
\usepackage{multirow}
\usepackage{epsfig}
\usepackage{soul}
\usepackage{colortbl}
\usepackage{pgf}
\usepackage{geometry}
\usepackage{times}
\usepackage{xspace}
\usepackage{epic}
\usepackage{pslatex}
\usepackage{graphicx}
\usepackage{caption}
\usepackage{subcaption}
\usepackage{hyperref}
\usepackage{tikz}

\usetikzlibrary{calc}
\usetikzlibrary{shapes.arrows}
\usetikzlibrary{positioning,fit,backgrounds}

\usepackage[textwidth=1.35\marginparwidth,
  backgroundcolor=yellow!50,
  linecolor=orange,
  textsize=scriptsize]{todonotes}

\definecolor{Red}{rgb}{1,0,0}
\definecolor{Gray}{rgb}{0.2,0.2,0.2}
\definecolor{Maroon}{rgb}{0.6,0.05,0.03}
\definecolor{Blue}{rgb}{0,0.7,0.9}
\definecolor{Green}{rgb}{0,.7,0}

\newcommand{\mads}{{\sc Mads}\xspace}
\newcommand{\pollster}{{\sc Pollster}\xspace}
\newcommand{\bimads}{{\sc BiMads}\xspace}
\newcommand{\multimads}{{\sc MultiMads}\xspace}
\newcommand{\dmultimads}{{\sc DMulti-Mads}\xspace}
\newcommand{\robustmads}{{\sc RobustMads}\xspace}
\newcommand{\psdmads}{{\sc PsdMads}\xspace}

\newcommand{\stomads}{{\sc StoMads}\xspace}

\newcommand{\nomad}{{\tt NOMAD}\xspace}

\newcommand{\sgtelib}{{\tt sgtelib}\xspace}

\newcommand{\orthomads}{{\sc OrthoMads}\xspace}

\newcommand{\stepstart}{{\texttt{Start}}\xspace}
\newcommand{\steprun}{{\texttt{Run}}\xspace}
\newcommand{\stepend}{{\texttt{End}}\xspace}
\newcommand{\worker}{{\sc Worker}\xspace}
\newcommand{\workers}{{\sc Workers}\xspace}

\newcommand{\compbiobjective}{{\texttt{BiObjective}}\xspace}
\newcommand{\compinitialization}{{\texttt{Initialization}}\xspace}
\newcommand{\compiteration}{{\texttt{Iteration}}\xspace}
\newcommand{\complh}{{\texttt{LH}}\xspace}
\newcommand{\compmegasearchpoll}{{\texttt{MegaSearchPoll}}\xspace}
\newcommand{\compneldermead}{{\texttt{NelderMead}}\xspace}
\newcommand{\compparallelspacedecomposition}{{\texttt{Parallel\-Space\-Decomposition}}\xspace}
\newcommand{\compmads}{{\texttt{Mads}}\xspace}
\newcommand{\comppoll}{{\texttt{Poll}}\xspace}
\newcommand{\compquadraticmodelsearch}{{\texttt{QuadraticModelSearch}}\xspace}
\newcommand{\compsgtelib}{{\texttt{SgtelibModelSearch}}\xspace}
\newcommand{\compsearch}{{\texttt{Search}}\xspace}

\newcommand{\compspeculativesearch}{{\texttt{SpeculativeSearch}}\xspace}
\newcommand{\compupdate}{{\texttt{Update}}\xspace}
\newcommand{\compvariableneighborhoodsearch}{{\texttt{VariableNeighborhoodSearch}}\xspace}

\newcommand\STOP{{\varepsilon_{\tt stop}}}
\newtheorem{algolocal}{Algorithm}[section]
\newenvironment{algo}[1]
	{	\begin{algolocal} {\em \bf #1} \hrulefill \end{algolocal}
		\begin{list}{}{} \item \vspace*{-2mm}			}
	{	\end{list} \vspace*{-3mm}\noindent \hrulefill \\			}

\newcommand\R{{\mathbb{R}}}
\newcommand\N{{\mathbb{N}}}
\newcommand\Z{{\mathbb{Z}}}

\begin{document}

\title{\nomad version~4: Nonlinear optimization with the \mads algorithm
\thanks{GERAD and Polytechnique Montr\'eal}
}

\author{
Charles Audet \thanks{Charles.Audet@gerad.ca}
\and
S\'ebastien {Le~Digabel} \thanks{Sebastien.Le.Digabel@gerad.ca}
\and
Viviane {Rochon~Montplaisir}\thanks{Viviane.Rochon.Montplaisir@gerad.ca}
\and
Christophe Tribes \thanks{Christophe.Tribes@gerad.ca}
}

%
%
%
%
%

\maketitle

\noindent
{\bf Abstract:}
\nomad is software for optimizing blackbox problems.
In continuous development since 2001,
it constantly evolved with the integration of new algorithmic features published in scientific publications.
These features are motivated by real applications encountered by industrial partners.
The latest major release of \nomad, version~3, dates from 2008. 
Minor releases are produced as new features are incorporated.
The present work describes \nomad~4, a complete redesign of the previous version,
 with a new architecture providing more flexible code, added functionalities and reusable code.
We introduce algorithmic components, which are building blocks for more complex algorithms, and can initiate other components, launch nested algorithms, or perform specialized tasks.
They facilitate the implementation of new ideas, including the \compmegasearchpoll component,
 warm and hot restarts, and a revised version of the \psdmads algorithm.
Another main improvement of \nomad~4 is the usage of parallelism, 
 to simultaneously compute multiple blackbox evaluations, 
 and to maximize usage of available cores.
Running different algorithms, tuning their parameters, and comparing their performance for optimization is simpler than before,
 while overall optimization performance is maintained between versions~3 and~4.
\nomad is freely available at
\href{https://www.gerad.ca/nomad}{\tt www.gerad.ca/nomad}
and the whole project is visible at
\href{https://github.com/bbopt/nomad}{\tt github.com/bbopt/nomad}. \\


\section{Introduction}

\nomad is software designed for the class of
{\em blackbox} optimization (BBO) problems~\cite{AuHa2017}. 
The term blackbox indicates that there is no information except the input and output. 
There is no analytic description of the objective 
    and/or constraint
functions, 
 there are no available derivatives, 
 possibly because they are not differentiable,
 and the functions may occasionally fail to return valid output and may require significant computational time to evaluate. 
This makes BBO problems difficult to solve, 
    in the sense that many optimization algorithms and heuristics cannot be applied. 
A typical BBO problem is a computer simulation of 
    an engineering or physical problem.

The development of \nomad was initiated in 2001 to implement direct search algorithms,
and major version~3 was released in 2008~\cite{Le09b}. The mesh adaptive direct search (\mads) algorithm~\cite{AuDe2006} is at the core of \nomad~3; it 
provides a flexible framework
and is supported by a rigorous hierarchical convergence analysis based on various degrees of smoothness of the functions defining the problem. 
Since its original release, minor releases of \nomad~3 have included several improvements and additions of algorithms to solve a variety of blackbox optimization problems efficiently. 

\nomad 
 and has proven its usefulness in scientific papers as well as in established companies.
Our own work includes contributions in 
    hydrology~\cite{AACW09a,SeCoAu2017,MiCaGuLeAuLe2014},
    pharmacology~\cite{SeCoAu2017},
    metamaterial design~\cite{AuLedDiSwMa2013},
    alloy design~\cite{GhLeAuCh2013,GhAuLeBeBaPe2012,GhRoLeAuPe2011},
    chemical engineering~\cite{AuBeCh2008a,HaBeAuKo03a}
    and bioinformatics~\cite{GeHiLedAuTerScha2014}.
Many other researchers use \nomad in a variety of fields.
In astrophysics for example, \nomad is used 
    for black hole observation~\cite{MTSMEBKF2019},
    for tracking the interstellar object 1I/{\textquoteleft}Oumuamua~\cite{HiHeEu2020},
    for kinematics analysis of galaxies~\cite{RNFSSSRGMR2019},
    and for gravitational wave detection, in a paper with more than 700 co-authors~\cite{PhysRevD.87.042001}.
Hundreds of applications are reported in the surveys~\cite{Audet2014a,AlAuGhKoLed2020,GhHaBeRoChPeBaLe2014},
    including many on energy, engineering design, and materials science.
    
In retrospect, the main development avenues of the \nomad
software and the \mads algorithm may be classified into three categories:
\begin{itemize}
    \item {\bf Algorithmic improvement.}
    The \mads algorithm was modified  to reduce the number of evaluations 
    through constraint handling techniques~\cite{AuDe09a,AuLedPey2014}, 
    by dynamically scaling and exploiting the specificity of variables~\cite{AuLe2012,AuLedTr2014,AuLeDTr2018},
    and by the improved integration of surrogates~\cite{AuCM2019,AuKoLedTa2016,sgtelib,TaAuKoLed2016,TaLeDKo2014}.
    
    \item {\bf Sub and super-algorithms.}
    \mads may call other optimization algorithms during its execution.  
    Sub-algorithms, when used under adequate conditions, may produce good candidate points for evaluation. Using the right candidates has a strong influence on the performance of the software. 
Useful algorithms and techniques are proposed in~\cite{AbAuDeLe09,AuIaLeDTr2014,AmAuCoLed2016,AuBeLe08,AuTr2018}. Sub-algorithms may also be tailored to exploit surrogate functions~\cite{AuKoLedTa2016,TaAuKoLed2016,AuCM2019}.
    Conversely, the \mads algorithm can be used as part of a broader direct search super-algorithm. For example, \bimads and \multimads~\cite{AuSaZg2008a,AuSaZg2010a} solve multiobjective optimization problems by running several \mads instances while managing the progress in obtaining a detailed Pareto front.
    \robustmads~\cite{AudIhaLedTrib2016} interrupts sequences of \mads runs by redefining the objective function to take into account noisy values.
    \psdmads~\cite{AuDeLe07} divides a large problem into problems of smaller dimension and launches instances of \mads in a parallel environment.
    
    \item {\bf Performance and parallelism.} 
    A major effort was placed into reducing the wall clock time to obtain good solutions. 
    The opportunistic strategy for evaluating points
    combined with ordering points to promote the most promising, ensure faster convergence~\cite{MScLASMC}. 
    Quadratic models approximate the problem to rapidly find better points~\cite{CoLed2011}. 
    Methods were developed to span a limited number of directions while maintaining the convergence proof, again to limit the number of blackbox evaluations~\cite{AudIhaLedTrib2016}.
    Conversely, to maximize core utilization during optimization, subspace exploration strategies~\cite{AuDeLe07,AdAuBeYa2014} as well as parallel strategies in \nomad~3~\cite{Le09b} were developed. However, these strategies are not well adapted to fully utilize the new abundance of computing resources, with some computer clusters counting cores in thousands.
\end{itemize}

Over the years, it became increasingly difficult to maintain and enhance the functionalities of \nomad~3. 
Recent algorithmic developments required modifications in many portions of the software.
The complex interactions between algorithms and sub-algorithms were not sufficiently anticipated. 
It was therefore decided to completely redesign the software.

The main goal of this new version of \nomad remains to solve efficiently a variety of constrained and unconstrained BBO problems.
In \nomad~4, the \mads algorithm as well as other algorithms deemed useful in \nomad~3 have been re-implemented 
    using primitive {\em algorithmic components}.
    which are building blocks for more complex algorithms, and interfaces adapted from the experience gained during the development of \nomad~3. 
This approach promotes software maintainability, as components may be reused when adding new algorithms.
This is an important requisite of this new version of the software even though it forces the rewriting of most of the source code.
The requisite to efficiently use a large number of available cores also had a strong impact on the architecture when redesigning the software.
Finally, the ability to tune algorithmic parameters (which control the algorithmic components) and to compare different algorithms is also an important requisite of the development. The optimization performance must be maintained between the versions. 

This paper describes the design of \nomad~4 to achieve this goal with the updated requisites. Sections~\ref{sec:mads}~and~\ref{sec:algocomp} present the \mads algorithm and other algorithmic components re-implemented from \nomad~3. The strategies for parallel blackbox evaluations in \nomad~4 are presented in Section~\ref{sec:parallel}.
The software architecture and development is presented in Section~\ref{sec:software:architecture}.
%
New algorithmic developments in \nomad must be assessed in terms of optimization performance and compared with other blackbox optimizers on a large variety of problems:
    Section~\ref{sec:results} compares the performance of the \nomad~3 and~4 versions and illustrates the gains produced by the use of multiple cores in \nomad~4.
Finally, Section~\ref{sec:conclusion} discusses future developments.

\section[The MADS algorithm]{The \mads algorithm}
\label{sec:mads}
\nomad solves optimization problems of the form
	$$\displaystyle \min_x \left\{f(x) ~:~ x \in \Omega \right\}$$
where $f: \R^n \rightarrow  \R \cup \{\infty \}$ is the objective function and $\Omega \subseteq \R^n$ is the feasible region.
Allowing the objective function to take the value $\infty$ is useful to exclude trial 
    points for which the evaluation failed to return valid output,
 for example, when the blackbox crashes or returns an error message.
The original \mads paper~\cite{AuDe2006}
    handled the constraint set by minimizing the unconstrained 
    extreme barrier function
	$f_\Omega : \R^n \rightarrow \R \cup \{\infty\}$ defined as
\begin{eqnarray*}
	f_\Omega(x) &:= &\left\{ \begin{array}{lll}
    f(x) &~\quad~ &\mbox{if}~x \in \Omega, \\
    \infty && \mbox{if}~x \notin \Omega.
    \end{array}\right.
\end{eqnarray*}
Later, 
    the progressive barrier~\cite{AuDe09a} approach was proposed
    to exploit the amount by which constraints are violated.
The optimization problem is redefined as
    \begin{equation}
             \min_{x \in \mathcal X  \subseteq \R^n} \left\{  f(x) ~:~  c(x) \leq 0 \right\},
       \label{pb-genctr}
    \end{equation}
where $f: \mathcal X \subseteq \R^n \rightarrow  \R \cup \{\infty \}$ and
      $c: \mathcal X \subseteq \R^n \rightarrow (\R \cup \{\infty \})^m$
      are functions with
      $c = \left(c_1, c_2, \dots, c_m \right)$, 
      and $\mathcal X$ is some subset of $\R^n$.
Again, the entire feasible region is denoted by
    $\Omega = \left\{ x \in \mathcal X ~:~ c(x) \leq 0 \right\}$.
The set $\mathcal X$ is frequently taken as being $\R^n$,
	the space of continuous variables,
	or as the set of nonnegative variables $\R^n_+$.

Each iteration of the \mads algorithm explores the space of variables through a global exploration called ``search step'', and a local exploration called ``poll step''.
Both these steps generate {\em trial points}, which are candidates for evaluation, on a discretization of $\mathcal X$ called the {\em mesh}.
At iteration $k$, 
    let $x^k$ denote the current best-known solution.
The mesh is defined as 
$M^k  :=   \left\{ x^k + \delta^k D y \ : \ y \in \N^p \right\} \ \subset \ \R^n$
where $\delta^k \in \R$ is the {\em mesh size parameter} and
 $D \in \R^{n \times p }$ is a positive spanning set of $p$ directions that satisfies  specific requirements.
The simplest possible set $D$ is the union of all positive and negative coordinate directions,
 and is the one implemented in \nomad, hence the following redefinition of the mesh at iteration $k$:
	\begin{eqnarray*}
		M^k & := &  \left\{ x^k + \delta^k y \ : \ y \in \Z^n \right\} \  \subset \ \R^n.
	\end{eqnarray*}

The search step is flexible, and allows the user to explore any finite number of mesh points in the set named $S^k$.
\nomad~4 proposes a one-point rudimentary line search in the direction of the previous success~\cite{AuDe2006},
    a Nelder~Mead inspired search step~\cite{AuTr2018}
    and a search based on the minimization of a quadratic model~\cite{CoLed2011,AmAuCoLed2016}.
Additional search strategies include basic Latin hypercube sampling~\cite{McCoBe79a}
    and others based on advanced statistical
    surrogates~\cite{AuKoLedTa2016,TaAuKoLed2016,TaLeDKo2014}. 
The {\em Variable Neighbourhood Search} search step~\cite{MlHa97a,AuBeLe08} available in \nomad~3 will be added in the future.
The user may also integrate their own search strategy.

The poll step follows more rigid rules than the search step.
Poll points are confined to a so-called frame around $x^k$
    whose dimension is set by a frame size parameter $\Delta^k$
    which is always greater than or equal to the mesh size parameter $\delta^k$.
The key elements of the poll step is that the poll set $P^k$ must lie within the frame,
    and $\left\{x - x^k : x \in P^k \right\}$ must be a positive spanning set for $\R^n$.
For \mads, as $k$ goes to infinity, the union of these normalized directions 
    becomes dense in the unit sphere.
\nomad~3 includes many examples of poll steps, 
    including coordinate search~\cite{FeMe1952},
    generalized pattern search~\cite{Torc97a},
    \orthomads with $2n$~\cite{AbAuDeLe09} and $n+1$~\cite{AuIaLeDTr2014} directions. 
Of these, \nomad~4 currently implements \orthomads with $2n$ directions;
 some simple direction strategies are also included;
 other types of poll steps will be included in the future. 

Algorithm~\ref{Algo-MADS algorithm} shows \mads with the extreme barrier to handle constraints. It is close to the one presented in the textbook~\cite{AuHa2017}.
\begin{algo}{
	\label{Algo-MADS algorithm}
	\sf 	Mesh adaptive direct search (\mads)}
Given $f : \R^n \rightarrow \R \cup \{\infty\}$, starting point $x^0 \in \Omega$, and barrier function $f_\Omega(x)$\\
{\sf 0. Initialization} \\
\hspace*{6mm}\begin{tabular}[t]{|lll}
		$\Delta^0 \in (0, \infty)$ 		&	& initial frame size parameter\\
		$\tau \in (0, 1),$ with $\tau$ rational 		&	& mesh-size adjustment parameter\\
		$\STOP \in [0, \infty)$ 		&	& stopping tolerance \\
		$k \leftarrow 0$				&	& iteration counter
	\end{tabular}\\ \\
{\sf 1. Parameter Update}\\
\hspace*{5mm} \begin{tabular}[t]{|l}
	set the mesh size parameter to $\delta^{k}  =  \min\left\{ \Delta^k, (\Delta^k)^2 \right\}$
\end{tabular}\\ \\
{\sf 2. Search}\\
\hspace*{5mm} \begin{tabular}[t]{|l}
	if $f_\Omega(t) < f_\Omega(x^k)$ for some $t$ in a finite subset $S^k$ of the mesh $M^k$ \\
	\qquad \begin{tabular}[t]{l}
 	 	set $x^{k+1} \leftarrow t$ and $\Delta^{k+1} \leftarrow \tau^{-1} \Delta^k$ and go to {\sf 4}\\
         \end{tabular} \\
         otherwise go to {\sf 3}\\
\end{tabular} \\ \\
{\sf 3. Poll}\\
\hspace*{5mm} \begin{tabular}[t]{|l}
    let $ P^k$ be a poll set constructed around $x^k$ using a positive spanning set\\
	if $f_\Omega(t) < f_\Omega(x^k)$ for some $t \in P^k$  \\
	\qquad \begin{tabular}[t]{l}
 	 	set $x^{k+1} \leftarrow t$ and $\Delta^{k+1} \leftarrow \tau^{-1} \Delta^k$ \\
         \end{tabular} \\
         otherwise \\
	\qquad \begin{tabular}[t]{l}
		set $x^{k+1} \leftarrow x^k$ and $\Delta^{k+1} \leftarrow \tau \Delta^k$ \\
	\end{tabular}
\end{tabular} \\ \\
{\sf 4. Termination}  \\
\hspace*{5mm} \begin{tabular}[t]{|l}
	if $\Delta^{k+1} \geq \STOP$ \\
\hspace*{5mm} \begin{tabular}[t]{l}	
	increment $k \leftarrow k+1$ and go to {\sf 1}
\end{tabular} \\	
otherwise stop and return $x^* \leftarrow x^{k+1}$
\end{tabular} \\
\end{algo}


The \mads algorithm with the progressive barrier to handle constraints is slightly more complex. The rules for accepting a new incumbent solution $x^{k+1}$ are based on both the objective function value and a second function that aggregates the constraint violations.
The rules also depend on whether $x^{k+1}$ is feasible or not.
A high-level description is found in Chapter~12 of~\cite{AuHa2017}
	and the detailed presentation appears in~\cite{AuDe09a}.

\section[Algorithmic components of \nomad]{Algorithmic components of \mads}
\label{sec:algocomp}

This section describes how the main elements of the \mads algorithm are encoded in \nomad~4.
The pseudo-code from Algorithm~\ref{Algo-MADS algorithm} offers flexibility and may be coded in different ways. 
The \mads algorithm takes as input an initial point $x^0$
    which is not required to be within the feasible region, 
    a set of algorithmic parameters $P$, 
    an objective function 
    $f$ and some constraint functions $c$.
These functions are provided through an executable code considered as a blackbox.
\mads provides the best incumbent solution $x^*$ found for this problem.

\subsection{Terminology}
\label{sec:algocomp:terminology}

When running an optimization program from given inputs, the code instructions are executed following a specific logic. The present section defines the terminology to describe the execution of \mads.

A {\em task} is defined as a group of code instructions adapted for human understanding. 
In our context, tasks can be generic or specialized, and they can be broken down into smaller tasks. 
Grouping tasks into algorithmic components allows to visualize the structure of the optimization program and its unfolding. 
 The name and purpose of a component come from the algorithm and the tasks that it performs.
The use of generic components and tasks during the design of an optimization program favors maintainability, which is a requisite for the revised version of \nomad.
Hence, the components can be readily reused in different optimization programs, and can be upgraded when new features are introduced.
%

The program consists of a hierarchy of components, which are executed one at a time, depth-first.
The components are named after the algorithms they implement, for example, \complh for the Latin hypercube sampling search step,
or \compmads for the \mads algorithm.
The execution of a component unfolds by performing the 
	{\em generic tasks} named \stepstart, \steprun and \stepend, in that order. 
{\em Specialized tasks} are tasks other than the generic tasks, with their name describing their purpose.
The \stepstart task may initiate another component, or perform a list of specialized tasks (from top to bottom in the figures below). 
The \stepend task may call specialized tasks required for the closure of the component.
\stepstart and \stepend tasks may also be void.
The \steprun task may initiate another component, or combine specialized tasks with iterating or testing.
Nested components and tasks model the execution of the entire program.

\subsection{Execution model for algorithms}
\label{sec:algocomp:execution}

\begin{figure}[htb!]
    \centering
     \includegraphics[width=0.8\textwidth]{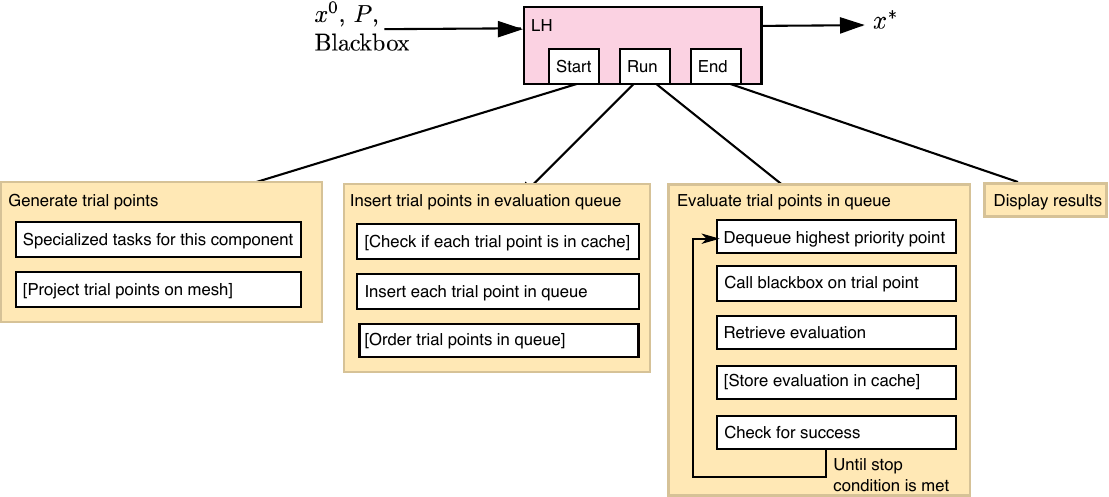}
    \caption{The \complh component (Latin hypercube sampling algorithm) with its connected tasks.
    Tasks within brackets are optional.}
    \label{fig:lhAlgo}
\end{figure}

\begin{figure}[htb!]
    \centering
    \includegraphics[width=0.8\textwidth]{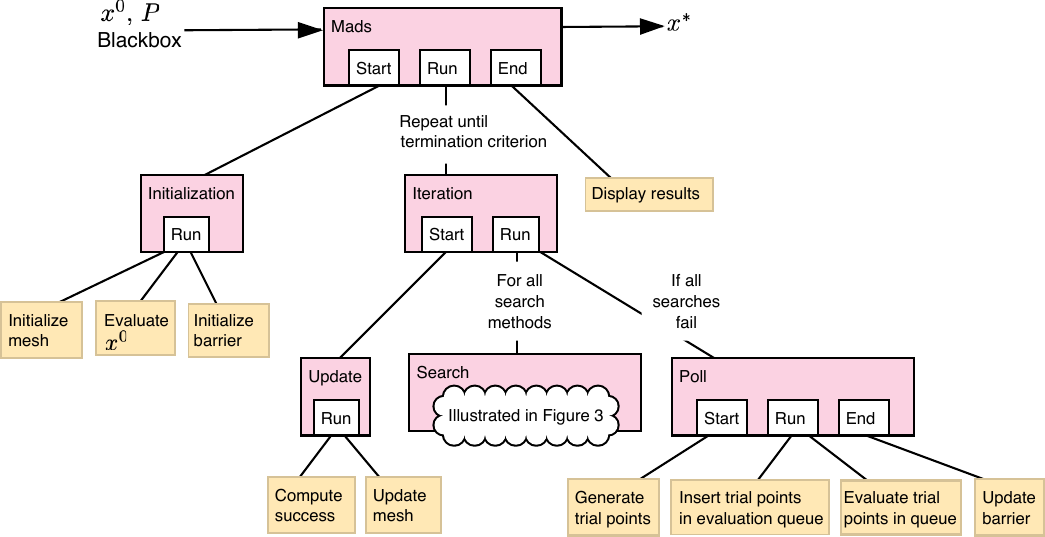}
    \caption{The \compmads component with its nested tasks and algorithmic components.}
    \label{fig:madsalgocomp}
\end{figure}

The execution model is illustrated in Figure~\ref{fig:lhAlgo}
 with the Latin hypercube sampling algorithm, and in Figure~\ref{fig:madsalgocomp}
 with the more complex \mads algorithm, which requires connecting several nested components and specialized tasks.
The \complh component presented in Figure~\ref{fig:lhAlgo}  sequentially performs four groups of specialized tasks: 
	{\tt Generate~trial~points} (called by \stepstart),
	{\tt Insert~trials~point~in~evaluation~queue} (called by \steprun), 	      {\tt Eva\-lu\-ate~trial~points~in~queue} (called by \steprun),
    and {\tt Display~results} (called by \stepend).
The \complh component can be called during the \mads search step (see Section~\ref{sec:mads}), or as a standalone optimization program. 
In both cases, all trial points are generated during the \stepstart task of the component. 

The \complh trial point generation details are found in~\cite{Tang93a}. 
Some tasks specialized for the generation and evaluation of trial points are presented in detail as they are common to all blackbox optimization programs. 
Some specialized tasks are optional, or are only necessary in some cases;
 they are written within square brackets.
For example, trial points are required to be located on the mesh when the \complh component is called during the \mads search step. 
In this case, trial points must be projected on the mesh prior to the actual evaluation. 
However, when the \complh component is called as a standalone algorithm, there is no mesh and therefore no projection is required.

When enabled, the {\em cache} contains the set of previously evaluated points, and the {\em incumbent} point is the best solution found yet. 
Each trial point is looked up in the cache. 
If it is not found, then it is added to the {\em evaluation queue} of points to be evaluated by the blackbox. 
Evaluations are run, possibly in parallel when multiple cores are available.
If a trial point evaluation is better than the incumbent point evaluation, it is possible to skip the points remaining in the queue
and to save the cost of evaluating them: this strategy is called {\em opportunism}.
There is a direct correspondence between some statements of Algorithm~\ref{Algo-MADS algorithm} and the tasks/sub-tasks. 
The comparison between $f_\Omega(t)$ and $f_\Omega(x^k)$
 seen in the search and poll steps corresponds to the tasks
 {\tt Call~blackbox~on~trial~point},
 {\tt Retrieve~evaluation}, and {\tt Check~for~success}.
%

The \stepstart task of the \compmads component executes the \compinitialization component referring to Step~0 in Algorithm~\ref{Algo-MADS algorithm}. 
As illustrated in the left part of Figure~\ref{fig:madsalgocomp},
 the \steprun task of the \compinitialization component first performs the mesh initialization and
 then conducts the provided initial point evaluation.
The task {\tt Evaluate~$x^0$} in  Figure~\ref{fig:madsalgocomp}
is identical to the previously described tasks
{\tt Insert~trials~point~in~eval\-uation~queue}
 and {\tt Evaluate~trial~points~in~queue}, with the trial point set to the initial point $x^0$. 
The \steprun task of the \compinitialization component
    concludes with the initialization of the progressive barrier parameters to handle the constraints.
The \steprun task of the \compmads component then repeatedly executes \compiteration components until a termination criterion is met. 
Each {\em \compiteration} component involves respectively an \compupdate, multiple \compsearch, and a \comppoll components. 
Each \compsearch component can initiate a sequence of nested algorithmic components. 
The \mads algorithm offers the flexibility to use any type of algorithm during the search step, 
    as long as a finite number of points is generated,
    that they are projected on the current mesh,
    and that the evaluation budget from the set of algorithmic parameters $P$ is not exceeded.

A component can call another instance of itself, directly or indirectly, as long as it is ensured that no infinite recursion is induced.
The nested components of the \compsearch component depend on which search step is performed. 
An example of search step based on quadratic models is described in 
 Section~\ref{sec:algocomp:combining} and Figure~\ref{fig:quadalgocomp}. 
The \comppoll and \compsearch components perform the tasks
 {\tt Generate~trial~points} (using different strategies),
 {\tt Insert~trial~points~in~evaluation~queue},
 {\tt Evaluate~trial~points~in~queue} and
 {\tt Update~barrier}.
While it is deployed, the program alternates generation and evaluation of trial points, 
 which can be a limiting factor to the number of parallel evaluations. 
Section~\ref{sec:parallel} presents a different way to deploy the execution of  \mads to exploit parallel blackbox evaluations.

\subsection{Combining algorithmic components}
\label{sec:algocomp:combining}
\begin{figure}[htb!]
    \centering
    \includegraphics[width=0.6\textwidth]{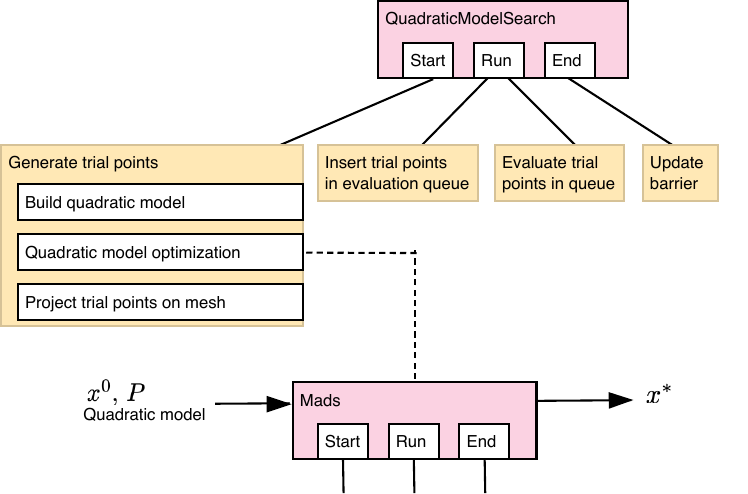}
    \caption{The \compquadraticmodelsearch component with its nested tasks and components. The nested components and tasks of the \compmads component are not presented.}
    \label{fig:quadalgocomp}
\end{figure}

Algorithmic components that represent algorithms can be run standalone. 
For instance, using \nomad, it is possible to run the \compneldermead algorithm to optimize a problem, 
 and the solution may be compared to the solution found using other algorithms such as \mads.
Algorithmic components are building blocks that may
    be combined and connected together to produce new algorithms.
The \complh component presented in Figure~\ref{fig:lhAlgo} can be used as a sub-algorithm by the \compsearch component of \mads
    to generate trial points, as long as the mesh projection is performed.
For the same purpose, we developed a \compquadraticmodelsearch~\cite{CoLed2011} component, in which
    previous blackbox evaluations are used to construct a quadratic model.
This model is used as a surrogate problem and is optimized to provide new trial points; this 
    optimization is performed by a new instance of \mads, with the quadratic model search step disabled to avoid infinite recursion.
 Therefore, we have a \compmads component which, through a \compquadraticmodelsearch component and the task
    {\tt Quadratic~model~optimization}, starts the execution of another instance of 
    the \compmads component. 
This is illustrated in Figure~\ref{fig:quadalgocomp}, with a dashed line connecting the two items. 

Various components can be used during a search with some control other the evaluation budget or the number of iterations. 
A \compspeculativesearch~\cite{AuDe2006} component generates trial points by using the direction of last success, which is speculated as a possible direction of improvement, starting from the current incumbent solution. 
We also developed a \compneldermead component~\cite{AuTr2018} based on the Nelder-Mead algorithm to iteratively generate and evaluate trial points. 
Version~3 of \nomad has a \compvariableneighborhoodsearch
component~\cite{AuBeLe08,MlHa97a,HaMl01a} to attempt escaping local solutions and will be integrated in \nomad~4 in the future. 
%

The \compmads component, with all its nested components and tasks, can itself be used within a super-algorithm that does not necessarily rely on a mesh: see for example~\cite{LiTr2017} where \nomad is hybridized with a mesh-free linesearch method.
It may also be executed repeatedly to solve a biobjective optimization problem through a series of single-objective formulations~\cite{AuSaZg2008a}. 
Version~3 of \nomad already has the corresponding \compbiobjective component that will be integrated in \nomad~4 in the future. 

Another example of the \compmads component being used as part of a super-algorithm is in \compparallelspacedecomposition (\psdmads~\cite{AuDeLe07}), where large problems are solved using an asynchronous parallel algorithm in which the parallel processes are launched on subproblems over subsets of variables.
A version of \psdmads is implemented in \nomad~4 using available algorithmic components. It is described in Section~\ref{sec:parallel:psdmads}.

\section{Parallel blackbox evaluations}
\label{sec:parallel}

A typical user of \nomad with access to a specific computational capacity 
    would like to obtain the best possible solution for an optimization problem within a certain time limit. 
This implies that \nomad must efficiently exploit all available cores. 
In some cases, the blackbox evaluation itself runs in parallel, 
    using all cores, but that is not always the case. 
An assumption for the software development is that running a blackbox evaluation requires significantly 
 more computational time compared to the other algorithmic tasks. 
Therefore, \nomad must efficiently distribute the blackbox evaluations in parallel, on secondary threads, while all other tasks are executed on a single main thread.
Several strategies for such parallel blackbox evaluations are presented in this section.

\subsection{The evaluation queue}
\label{sec:parallel:queue}
Regardless of the strategy used to manage the parallel evaluations, 
	an {\em evaluation queue} is maintained to manage evaluations, in which the elements are trial points to be evaluated. 
It behaves as a priority queue, and is implemented as a sorted vector. 
When a trial point is generated, 
    it is inserted in the queue, provided that it was not previously evaluated.
The trial points may be ordered, 
    so that the most promising one are evaluated first. 
Sorting the points is important when opportunism is used,
because when the evaluation of a point leads to a new success, the remaining points in the queue are not evaluated.
Different ordering strategies are available to sort the trial points in the queue:
    based on the direction of the last success, 
    on the order in which trial points were generated,
    on the lexicographic order of their coordinates, or simply mixed randomly.
    
\subsection{Grouping evaluations}
\nomad provides the option to group trial points, with a given maximum group size, for evaluation.
With this strategy, users are in charge of managing the dispatching 
  of the groups of points,
  depending on the specifics of their blackbox and computers, in order to maximize core usage.
Nevertheless, grouping points for evaluation is not ideal 
    because there may be an insufficient number of trial points in the queue to fill a group to its maximum size, 
    resulting in unexploited cores.
Exploratory work on filling groups of poll sets appears in~\cite{MScGL} and will eventually be incorporated into \nomad.
 



\subsection{Parallel evaluations on multiple threads}
\label{sec:parallel:mainthread}
\nomad manages one or more {\em main threads} and, if additional cores are available, optional {\em secondary threads}.
In this Subsection we consider a single main thread. 
See Section~\ref{sec:parallel:psdmads} for a case where multiple main threads are used.
The main thread performs all algorithmic tasks, including some of the evaluations of trial points.
The secondary threads only execute evaluation tasks.
These tasks, independently of the algorithm, are grouped under the name ``Evaluate \textbf{one} trial point in queue'', and are the same as those listed in ``Evaluate trial point\textbf{s} in queue'' in Figure~\ref{fig:lhAlgo}. 
Each thread performs one evaluation task at a time.

    

\begin{figure}[htb!]
    \centering
    \includegraphics[width=0.5\textwidth]{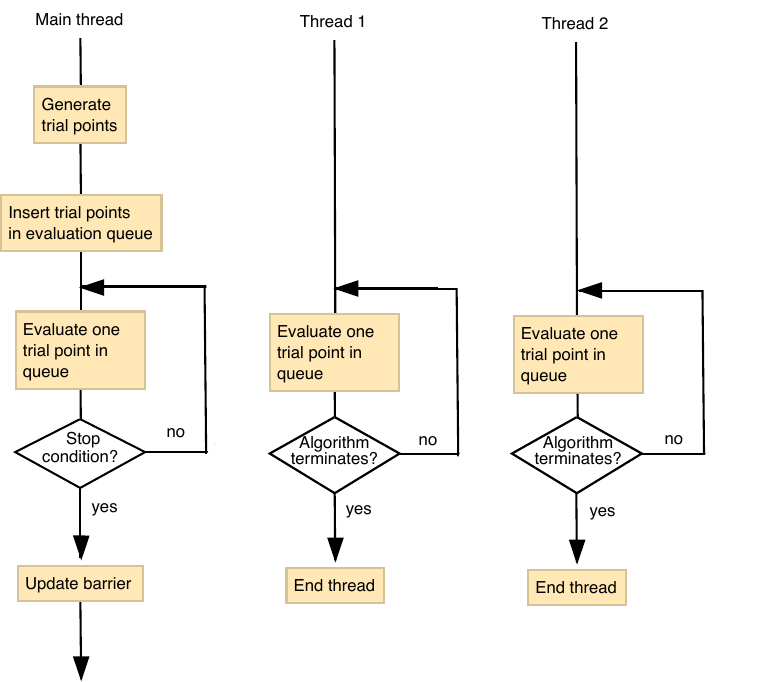}
    \caption{Workflow of the tasks for the evaluation queue using parallel threads.}
    \label{fig:allThreads}
\end{figure}

Figure~\ref{fig:allThreads} illustrates the tasks workflow performed in the main thread and two secondary threads. 
The specialized tasks for an algorithm are not presented. 
In the main thread, the iterative process of dequeuing and evaluating points terminates when there are no more points in the queue, or when some criteria is met 
    (budget of evaluations is filled, opportunism condition is met, etc.)
At this stage, the evaluation queue is cleared from any remaining unevaluated trial points,
the progressive barrier parameters are updated, 
    new trial points are generated, and the algorithm continues. 
In the meantime, the secondary threads keep on working and
    the iterative process of dequeuing and evaluating points goes on until the algorithm terminates in the main thread. 
All information relative to a trial point 
    (success, objective and constraint values) 
    that is evaluated in a secondary thread is made available to the other tasks of the algorithm, and this information is taken into account in the main thread for the continuation of the algorithm. 

The management of threads is currently done using {\tt OpenMP}. 
The user may provide the number of threads $n_t$ to efficiently access the computer cores;
    otherwise, {\tt OpenMP} computes the number of available threads.

\subsection{Grouping trial points generation}
\label{sec:parallel:megasearchpoll}

\begin{figure}[htb!]
    \centering
    \includegraphics[width=0.8\textwidth]{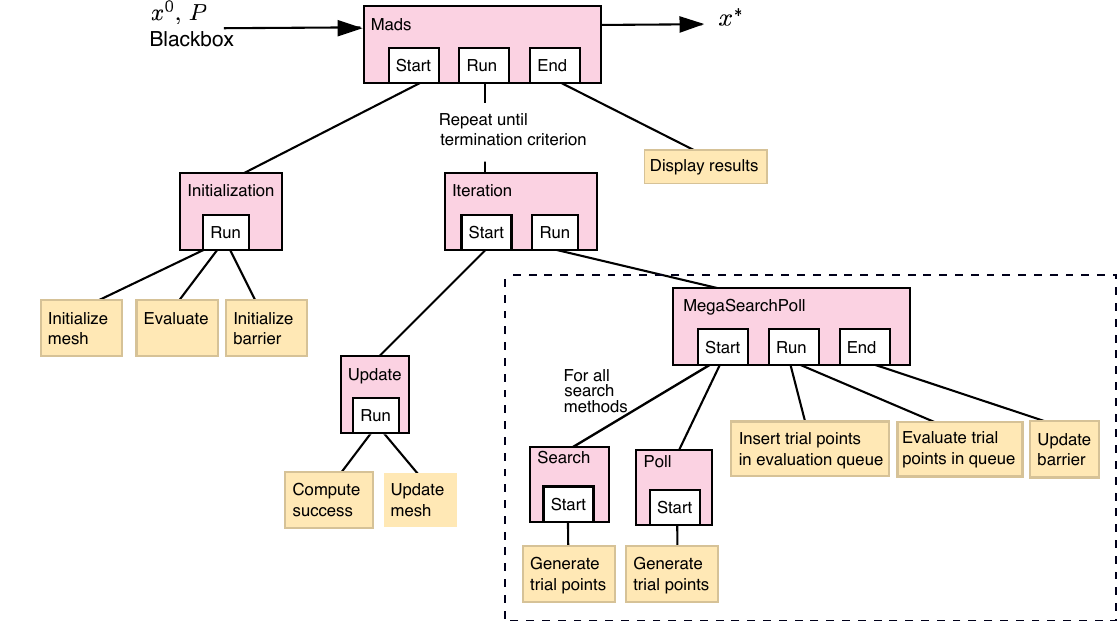}
    \caption{The \compmegasearchpoll version of \mads with its nested tasks and algorithmic components. A \compmegasearchpoll component is added to the \compmads component of Figure~\ref{fig:madsalgocomp}. The \compmegasearchpoll component and its dependant tasks and components are presented in a dashed box. The \steprun and \stepend tasks of the \compsearch and \comppoll are not called.}
    \label{fig:megasearchpoll}
\end{figure}

In every \compsearch and \comppoll components of the generic \mads algorithm described in
 Section~\ref{sec:mads}, evaluations are performed immediately after the trial points are generated. 
This approach generates few points, sometimes a single point, to be evaluated, 
 which makes it difficult to exploit multiple cores for evaluations and to use parallelism at its full capacity. 
We developed a new combination of \compsearch and \comppoll components, 
 called {\em \compmegasearchpoll} (see Figure~\ref{fig:megasearchpoll}).
It generates all the trial points for the \compsearch and \comppoll components,
  and only then the points are inserted in the queue and evaluated in parallel. 
This way, more points are evaluated at a given time. 
Additionally, search strategies like the speculative search were reworked to provide more points. 
Different strategies to enlarge the number of points generated by the poll step are also proposed and examined in~\cite{MScGL}. 

\subsection[\psdmads: \mads with parallel space decomposition]{\psdmads: \mads with parallel space decomposition}
\label{sec:parallel:psdmads}
\begin{figure}[htb!]
    \centering
    \includegraphics[width=0.7\textwidth]{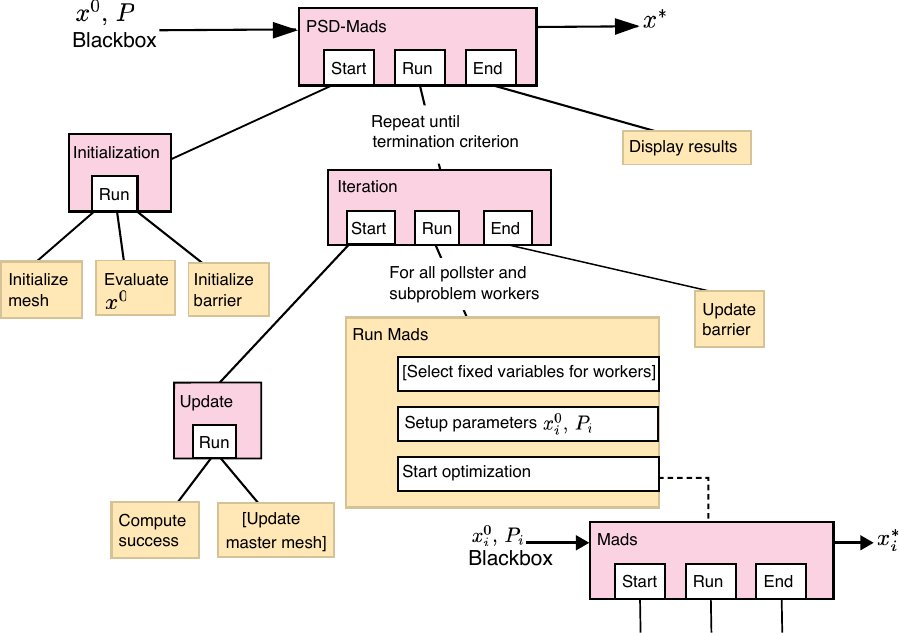}
    \caption{Description of \psdmads using algorithmic components. The \pollster and subproblem \workers are run in parallel. The nested components and tasks of \compmads are not presented.}
    \label{fig:psdmadsalgocomp}
\end{figure}

The sequential implementation of \mads is recommend from problems whose dimension is reasonably small ($n \leq 50$).
\psdmads was developed~\cite{AuDeLe07} to solve larger problems using space decomposition and parallelism.
In \psdmads, the problem is divided into random subspaces of dimension $n_s$, much smaller than $n$. 
Values of $n_s$ ranging from $2$ to $4$ are frequently used.  
Each subproblem is solved by launching a \mads 
  algorithmic component
  called \worker.
An additional \mads algorithmic component
    called \pollster is launched in the space of dimension $n$, 
    but evaluates a single point generated by a \comppoll component.
In practice, the \pollster rarely improves the incumbent solution but
its presence is necessary to ensure that the theoretical convergence results of \psdmads are satisfied. 
The \pollster and \workers are repeatedly 
 launched 
 within an \compiteration component. 
The mesh sizes of the \pollster and \workers are bounded by a master mesh that
is updated at every \compiteration.
These bounds, the \pollster's single evaluation, and the \worker's small dimension, allow for a fast resolution of each \mads.
Figure~\ref{fig:psdmadsalgocomp} presents the main algorithmic components involved in \psdmads.

\begin{figure}[htb!]
    \centering
    \includegraphics[width=0.6\textwidth]{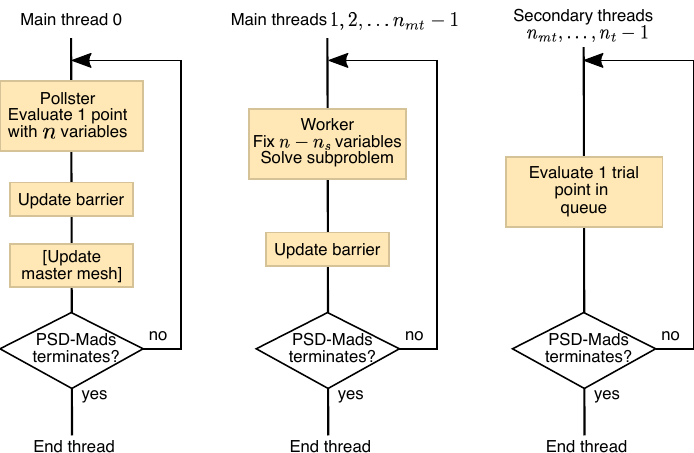}
    \caption{Thread usage in \psdmads. $n_s$ is the number of variables handled by each \worker.}
    \label{fig:psdmadsthreads}
\end{figure}

In the original \psdmads implementation described in~\cite{AuDeLe07},
the management of parallel processes is done using {\tt Message Passing Interface (MPI)}.
In the new implementation, {\tt OpenMP} manages the parallel execution of algorithmic components on main threads.
The total number of available threads is $n_t$.
The \pollster and \workers are run by \compmads components on $n_{mt}$ main threads, with $n_{mt} \leq n_t$ as illustrated by  Figure~\ref{fig:psdmadsthreads}. 
Main thread $0$ is used for the \pollster.
Main threads $1$ to $n_{mt}-1$ are used for the \workers.
Additional threads $n_{mt}$ to $n_t-1$ are secondary threads.    

The main threads are used for algorithmic components and evaluations, whereas the secondary threads are used only for evaluations.
The master mesh size update (enlarged or refined) depends on the success of one of the \workers or \pollster.
In the new implementation, a finer update control 
 delays the mesh size update until a prescribed minimum number of variables are explored by subproblems solved by \workers.


\section{Software architecture and development}
\label{sec:software:architecture}

This section reviews high level choices made for the different software components of \nomad.
Next, details relative to the processes and tools used while creating this software are mentioned.

\subsection{Software architecture}
The previous section illustrated how \nomad~4 and algorithmic components facilitate the construction of new algorithms, such as \psdmads.

\begin{figure}[htb!]
    \centering
    \includegraphics[width=0.7\textwidth]{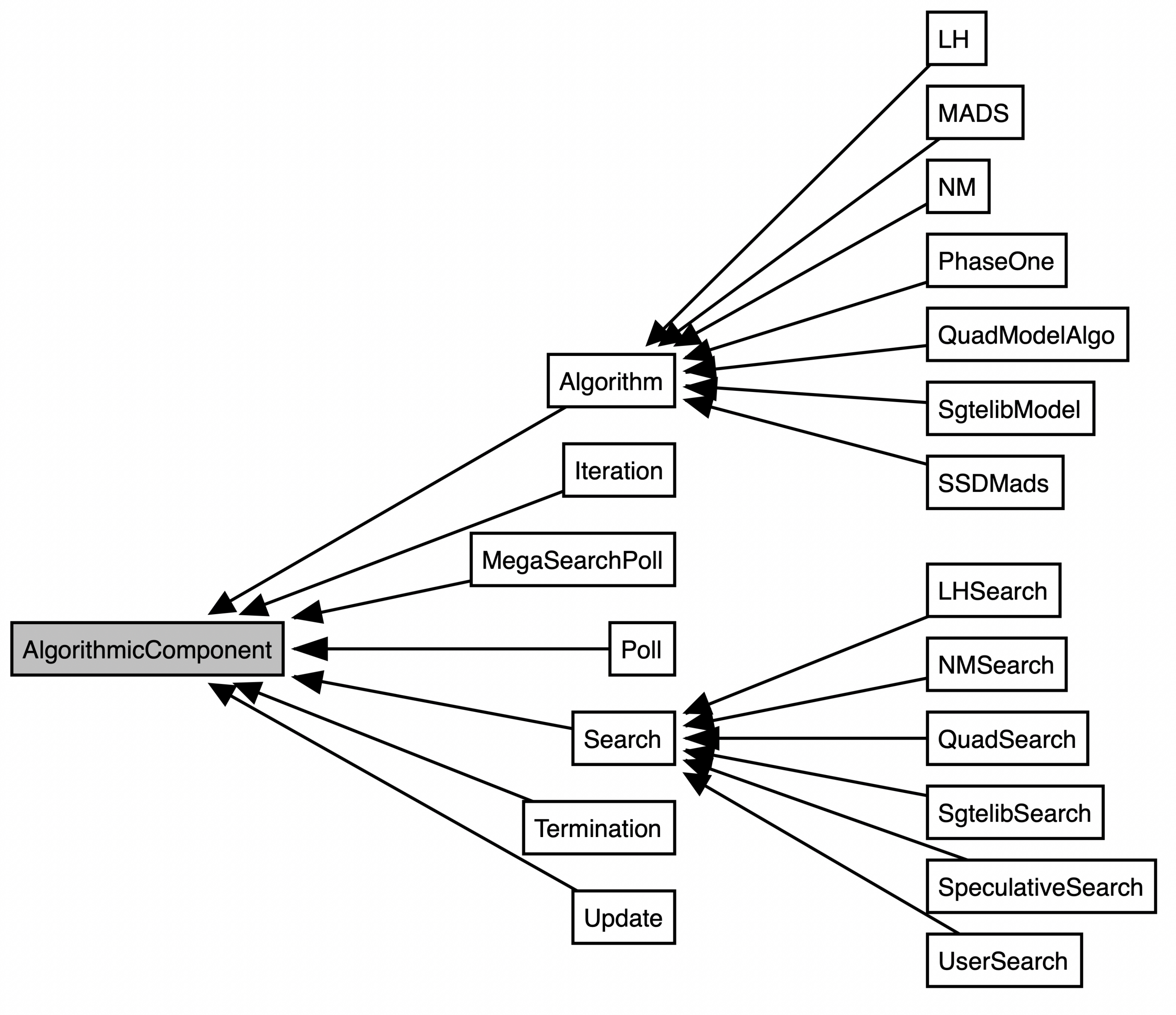}
    \caption{Inheritance graph for the main {\tt AlgorithmicComponent} classes.}
    \label{fig:stepinheritance}
\end{figure}


Since this new version is written from scratch, an effort was placed in the development of a modular and reusable architecture.
The algorithmic components described in Section~\ref{sec:algocomp} are building blocks for algorithms.
The evaluation queue presented in Section~\ref{sec:parallel} is designed to launch evaluations in parallel.
The code is implemented using object-oriented programming.
Algorithmic components and evaluation queue are translated in code as classes and objects.
Polymorphism is used:
 for instance, all algorithmic components are of the base class {\tt Step}; classes {\tt Algorithm} and {\tt SearchMethod} inherit from {\tt Step}; 
 class {\tt SpeculativeSearch} derives from {\tt SearchMethod}. 
Efforts are made to ensure that the code is clear and generic enough to be easily understood and modified, 
 for example by implementing a new search method.

The code is organized into libraries that provide a range of functionalities for programming algorithms to solve optimization problems:

\begin{itemize}
    \item
    {\tt Utils}: Math functionalities; parameter definition and checking;
    output management, including an output queue to correctly display information coming from different threads;
    file utilities, clock, and other utilities.
    \item
    {\tt Eval}: All that relates to the management of evaluations: Evaluation queue, evaluator, results of evaluations, and cache for points that have already been evaluated.
    \item
    {\tt Algos}: Algorithmic components and algorithms: \complh, \compmads, \compneldermead, \compquadraticmodelsearch, 
    \compparallelspacedecomposition, and \compsgtelib, where the \compsgtelib is a more general case of the \compquadraticmodelsearch algorithm.
    \item
    \sgtelib~\cite{sgtelib}, a library containing many surrogate models developed by B. Talgorn. 
\end{itemize}

\subsection{Software development}
The development process of \nomad~4 is inspired by the Agile software development values.
The team meets over daily scrums and biweekly group meetings. 
They work closely with students and business partners.
Features and issues are discussed and added timely to the code.
Code quality is verified through unit tests, for classes and methods, and through integration tests, for algorithmic functionality.
Performance profiles (presented in section~\ref{sec:results}) are processed regularly,
    comparing \nomad~4 with \nomad~3 or with previous versions of \nomad~4,
    to establish that development is going in the right direction to efficiently solve optimization problems.

Customer collaboration is key to development. 
For instance, one of our key users asked for {\em hot/warm restart}.
This new feature makes it possible to continue the solving process after it has started, without having to restart it
 from the beginning.
In the case of hot restart, the user interrupts the solver to change the value of a parameter. 
With warm restart, the user changes a parameter from a resolution that has already reached a termination condition. 
In both cases, the solving process is then continued from its current state.
This feature was discussed with the user, and added to \nomad~4. The user could test it promptly.

\nomad~4 is a standalone program coded in {\tt C++14} using {\tt OpenMP} when available,
on {\tt Linux} and {\tt macOS}.
A {\tt Windows} version will be available soon.
{\tt CMake} is used for compilation.
\href{https://github.com/google/googletest}{\tt Google Test} is used for unit tests.
Stable code is available and updated frequently at \href{https://github.com/bbopt/nomad}{\tt github.com/bbopt/nomad}.


\section{Computational results}
\label{sec:results}
The code of \nomad~4 differs significantly from that of \nomad~3;
    only a few base classes were preserved.
Comparing the performance of the two versions is crucial to validate 
    that algorithms have been correctly coded.

Tests are conducted using an in-house application called the Runner.
The Runner is designed to compare the 
    performance of different optimization software including
    different versions of \nomad and various algorithmic choices. 
The benchmark tests presented below
    include constrained and unconstrained analytical problems,
    engineering test problems from the literature,
    as well as tests involving parallelism.

Comparisons are made through
    data profiles~\cite{MoWi2009}.
The vertical axis shows the proportion of problem solved 
    within a prescribed tolerance of a parameter $\tau$
    and the horizontal axis measures the effort deployed by the compared methods in
    terms of groups of $n+1$ function evaluations.
A steep curve indicates that the corresponding method rapidly improves the solutions.
A method having its curve above the others performs better for the prevailing test conditions.
The optimization runs are conducted on a series of problems for a given evaluation budget. 
In the tests below, 
    each graph has two curves, one for \nomad~3 and the other for \nomad~4 with all default parameters, except that the direction type is set to ORTHO~2N  (\nomad~3 has a default called ORTHO~N+1 QUAD direction type that is not yet implemented in \nomad~4) 
    and the ordering of the points before evaluation is set to the last direction of success 
    (\nomad~3 has a default ordering based on quadratic models that is not yet implemented in \nomad~4).

\subsection{Validation on analytical problems}

\begin{figure}[htb!]
    \centering
    \begin{subfigure}[b]{0.45\textwidth}
         \centering
         \includegraphics [height=5.8cm, ext=pdf]{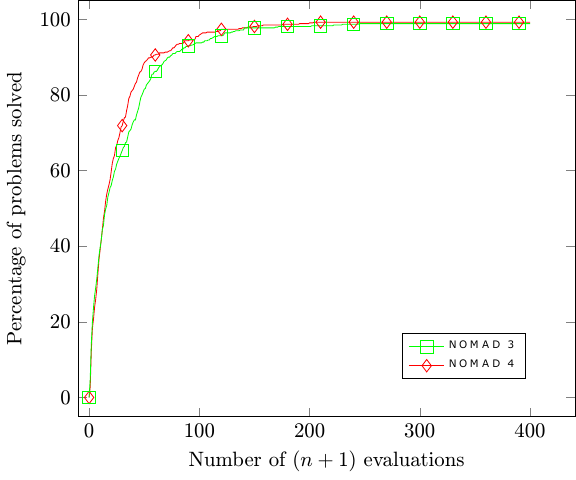}
         \caption{$\tau=10^{-2}$.}
         \label{fig:dp2_53pbsNoCons_MW_smooth}
    \end{subfigure} 
    \hfill
    \begin{subfigure}[b]{0.45\textwidth}
         \centering
         \includegraphics [height=5.8cm, ext=pdf]{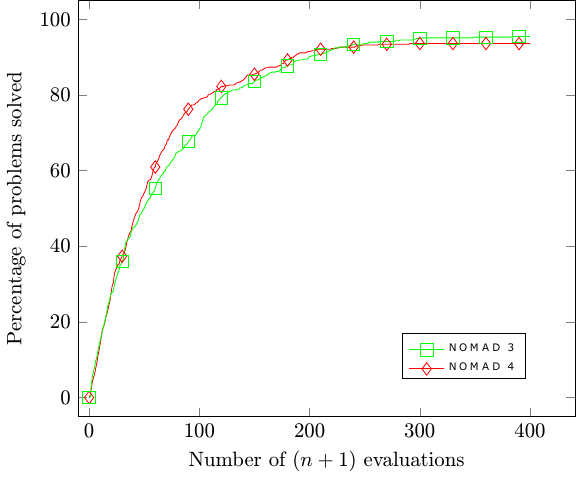}
         \caption{$\tau=10^{-4}$.}
        \label{fig:dp4_53pbsNoCons_MW_smooth}
    \end{subfigure}
    \caption{Data profiles obtained on 53 smooth unconstrained problems. The precision to detect if a problem is solved is set to $\tau=10^{-2}$ (left) and $\tau=10^{-4}$ (right).} 
    \label{fig:dataprofiles:MW}
\end{figure}

Figure~\ref{fig:dataprofiles:MW}
    compares \nomad~3 and~4 on a collection of 53 unconstrained smooth problems~\cite{MoWi2009}
    with a number of variables $n$ ranging
    from 2 to 12 
    and with a budget of $400(n+1)$ function evaluations.
Each version is launched 10 times with different random seeds.
The profiles on the left use a tolerance of $\tau = 10^{-2}$
    (the reader is invited to consult~\cite{MoWi2009}
        for the precise description of the role of $\tau$)
    and the ones on the right use a smaller tolerance of $\tau = 10^{-4}$.
In both cases, the two versions exhibit comparable performance,
    and there is no clear dominance of one over the other.

\begin{figure}[htb!]
    \centering
    \begin{subfigure}[b]{0.45\textwidth}
         \centering
         \includegraphics [height=5.8cm, ext=pdf]{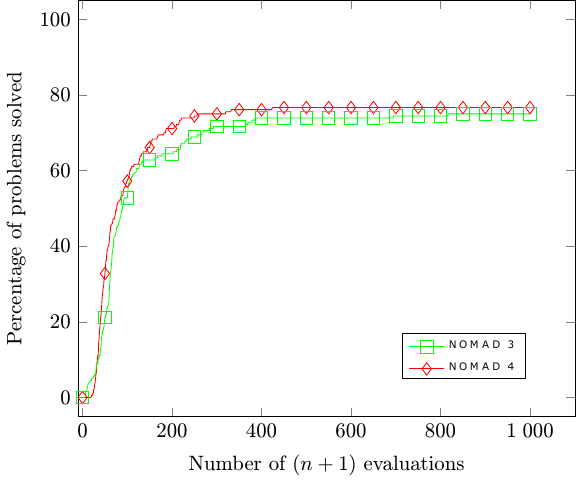}
         \caption{$\tau=10^{-2}$.}
         \label{fig:dp2_18pbsWithCons}
    \end{subfigure} 
    \hfill
    \begin{subfigure}[b]{0.45\textwidth}
         \centering
         \includegraphics [height=5.8cm, ext=pdf]{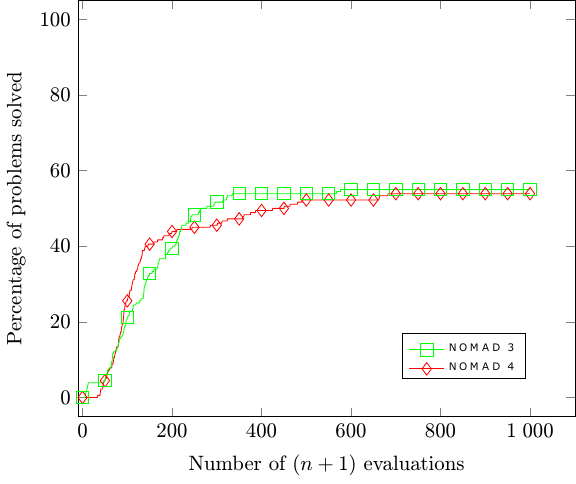}
         \caption{$\tau=10^{-4}$.}
        \label{fig:dp4_18pbsWithCons}
    \end{subfigure}
    \caption{Data profiles obtained on 18 problems with constraints. The precision to detect if a problem is solved is set to $\tau=10^{-2}$ (left) or $\tau=10^{-4}$ (right).} 
    \label{fig:dataprofiles:18PbsWithCons}
\end{figure}

Figure~\ref{fig:dataprofiles:18PbsWithCons}
    compares \nomad~3 and \nomad~4 on the collection of 18 constrained problems listed in Table~\ref{tab-pbs} with a budget of $1000(n+1)$ function evaluations.
The number of variables ($n$) varies from 2 to 20,
    the number of constraints ($m$) ranges
    from  1 to 15,
    and 13 of the problems have bounds on the variables.
Again, each version is launched with 10 random seeds.
The profiles on the left use a tolerance of $\tau = 10^{-2}$
    and the ones on the right use a smaller tolerance of $\tau = 10^{-4}$.
As with the unconstrained case, \nomad~3 and~4 have a similar performance.

\begin{table}[htb!]
\begin{center}
\renewcommand{\tabcolsep}{3pt}
\begin{footnotesize}
\begin{tabular}{|rrrrrc|}  
\hline
\# & Name & Source & $n$ & $m$ & Bnds \\ 
\hline
\hline
 1 &  CHENWANG\_F2 	&\cite{ChWa2010} & 8 & 6 & yes  \\ 
 2 &  CHENWANG\_F3 	&\cite{ChWa2010} & 10 & 8 & yes \\ 
 3 &  CRESCENT 		&\cite{AuDe09a} & 10 & 2 & no   \\ 
 4 &  DISK  		&\cite{AuDe09a}  & 10 & 1 & no  \\ 
 5 &  G210			&\cite{AuDeLe07}             & 10 & 2 & yes  \\ 
 6 &  G220			&\cite{AuDeLe07}                & 20 & 2  & yes \\ 
 7 &  HS19			&\cite{HoSc1981} & 2 & 2 & yes  \\ 
 8 &  HS83 			&\cite{HoSc1981} & 5 & 6 & yes  \\ 
 9 &  HS114 		&\cite{LuVl00}               & 9 & 6 & yes  \\ 
10 &  MAD6			&\cite{LuVl00}                & 5 & 7 & no  \\ 
11 &  MDO 			&\cite{TriDuTre04a}                & 10 & 10  & yes\\ 
12 &  MEZMONTES		&\cite{MezCoe05} & 2 & 2 & yes  \\ 
13 &  OPTENG\_RBF 	&\cite{KiArYa2011}       & 3 & 4 & yes  \\ 
14 &  PENTAGON		&\cite{LuVl00} & 6 & 15    & no   \\ 
15 &  SNAKE 		&\cite{AuDe09a} & 2 & 2 & no          \\ 
16 &  SPRING 		&\cite{RodRenWat98} 	& 3 & 4 & yes \\ 
17 &  TAOWANG\_F2 	&\cite{TaoWan08} & 7 & 4 & yes  \\ 
18 &  ZHAOWANG\_F5 	&\cite{ZhaWan2010b} & 13 & 9 & yes \\ 
\hline
\end{tabular}
\end{footnotesize}
\end{center}
\caption{Description of the set of 18 analytical problems with constraints.}
\label{tab-pbs}
\end{table}

\subsection[Tests on PSD-MADS]{Tests on \psdmads}
The implementations of \psdmads in \nomad~3 and in \nomad~4 are significantly different,
 principally because of the strategy used to perform parallel subproblem optimizations. 
We present a comparison of the two \psdmads implementations on two bound constrained variants of the Rosenbrock test problem~\cite{GoOrTo03}. 
The two variants, called SRosenbr50 and SRosenbr250, have $n=50$ and $n=250$ variables respectively, 
    with all starting point coordinates set to 0.5. 
The lower bounds are all set to -10 and the upper bounds are set to 10. 

\begin{figure}[htb!]
    \centering
    \begin{subfigure}[b]{0.45\textwidth}
         \centering
         \includegraphics [height=5.8cm, ext=pdf]{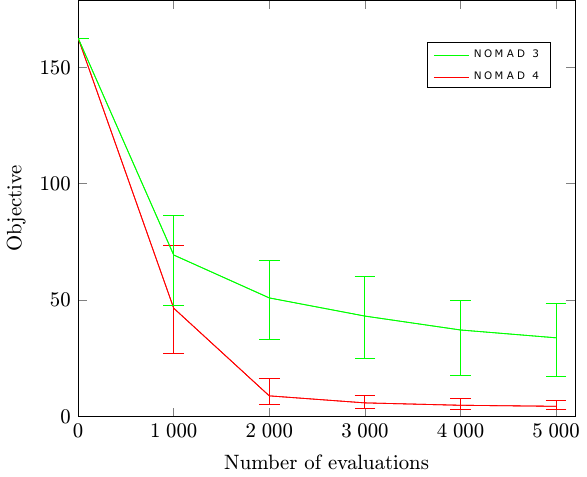}
         \caption{SRosenbr50.}
         \label{fig:srosenbr50}
    \end{subfigure} 
    \hfill
    \begin{subfigure}[b]{0.45\textwidth}
         \centering
         \includegraphics [height=5.8cm, ext=pdf]{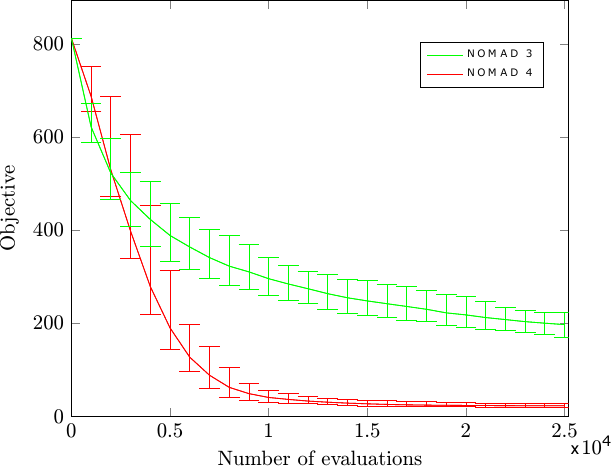}
         \caption{SRosenbr250.}
        \label{fig:srosenbr250}
    \end{subfigure}
    \caption{Convergence plot for 30 runs of \psdmads on 
    SRosenbr50 (left) and SRosenbr250 (right). The solid line is for the best average objective function value and the error bars represent the minimum and maximum objective function values after 1, 1000, 2000, \ldots, 100$n$ times the number of blackbox evaluations.}
    \label{fig:plotPSDMads}
\end{figure}

Because of the stochastic nature of \psdmads, both instances of the problem are solved 30 times each to perform a fair comparison. 
Figure~\ref{fig:plotPSDMads} 
    plots the average incumbent objective function value
    versus the number of function evaluations.
The plot also shows the best and worst objective function values.


The two implementations of \psdmads use 4 workers to perform subproblem optimizations with 2 
    randomly selected variables among the 60 available. 
\nomad~3 uses {\tt MPI} with 6 processes (1 process for the pollster, 3 processes for the regular workers,
    1 process for the cache server and 1 process for the master). 
To obtain a comparable task distribution, \nomad~4 uses only $n_t=4$ main threads for {\tt OpenMP} (no secondary thread) for the pollster and regular workers (see Figure~\ref{fig:psdmadsthreads}).

For both variants, the \psdmads version in \nomad~4
    outperforms that of \nomad~3.
The worst performance of the 30 \psdmads runs using \nomad~4 
    is significantly better than the best run of \nomad~3.

\subsection{Improvements in solution times}
The previous section compared the performance of \psdmads 
    in terms of number of function evaluations.
We next study the impact of using multiple cores on the overall computational time.
The Rosenbrock test problem is not adequate for such comparisons, 
    as it is evaluated nearly instantaneously.
We present results on the problem Solar~7~\cite{MScMLG}, which requires approximately 5~seconds for each evaluation.
This problem simulates the operation of a solar thermal power plant.
It has 7 variables, one of which is integer, 6 constraints, and variables are bounded.

\begin{figure}[htb!]
    \centering
    \begin{subfigure}[b]{0.45\textwidth}
         \centering
         \includegraphics [height=5.8cm, ext=pdf]{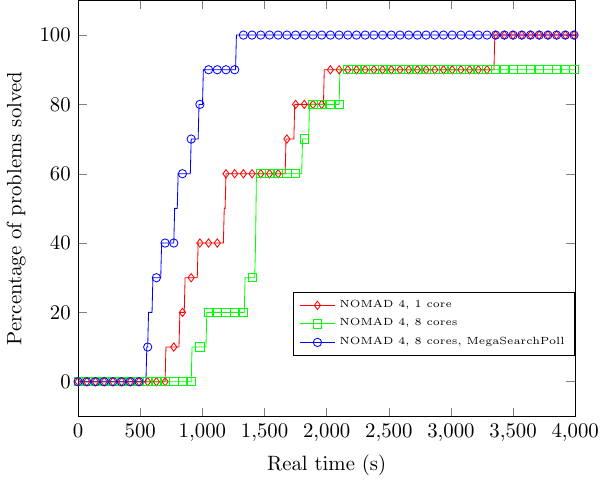}
         \caption{Data profiles with respect to time; $\tau = 10^{-2}$.}
         \label{fig:timeprofiles:tdp2}
    \end{subfigure} 
    \hfill
    \begin{subfigure}[b]{0.45\textwidth}
         \centering
         \includegraphics [height=5.8cm, ext=pdf]{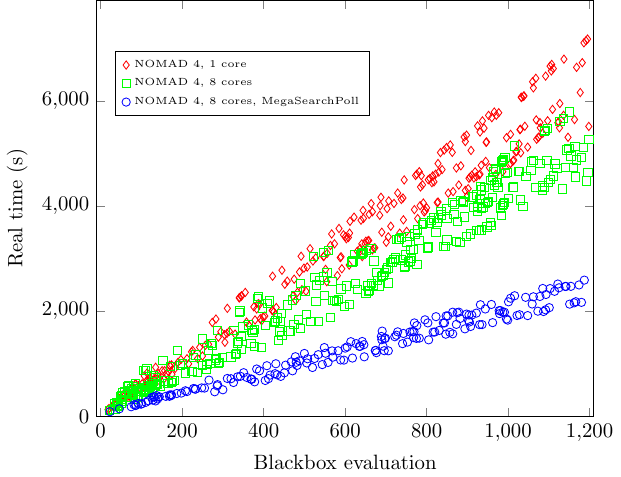}
         \caption{Speed-up.}
         \label{fig:timeprofiles:timeprofile}
    \end{subfigure}
    \caption{Effect of using multiple cores for \nomad~4, for constrained problem Solar~7,
    where a single evaluation takes about five seconds.} 
    \label{fig:timeprofiles}
\end{figure}

Figure~\ref{fig:timeprofiles:tdp2} shows 
    data profiles, for different parameter settings of \nomad~4,
    where the $x$-axis represents the wall-clock time in seconds
    rather than the number of function evaluations.
A method having its corresponding curve above the others performs faster than the others.
Figure~\ref{fig:timeprofiles:timeprofile} illustrates the speed-up,
    which plots the wall clock time in seconds 
    in terms of the number of function evaluations. 
Low values on the plot indicate better performance.
Three cases were tested on a machine containing 8~cores, by varying $n_t$, the number of threads used,
    which here is equal to the number of cores used.
In the first case, a single core is used ($n_t$ = 1).
In the second case, 8~cores are used ($n_t$ = 8).
The third case combines 8~cores with the \compmegasearchpoll component
    described in Section~\ref{sec:parallel:megasearchpoll}. 
Each case is launched ten times with different random seeds.

Without surprise,
 using multiple cores in parallel allows \nomad~4 to 
 generate solutions faster.
In addition, using the \compmegasearchpoll component
    improves the speed even more.
In summary, running \nomad~4 on Solar~7 
    using 8~cores as well as the \compmegasearchpoll component 
    leads to an overall computational time up to 3.3~times faster 
    than using a single core.

\section{Conclusion}
\label{sec:conclusion}
The \nomad blackbox optimization package has been completely redesigned. 
The new design defines easily interfaceable building blocks named algorithmic components, 
    for constructing elaborate algorithms.
This approach promotes software maintainability and modularity.
The architecture is strongly impacted by the requirement of efficiently using a large number of cores.

The new version's numerical performance is comparable to that of the previous version.
However,
 the code is easily accessible for students and developers.
The modularity of the algorithmic components makes the code flexible, reusable and allows for easy development of new algorithms.
Because the architecture is designed with parallelism in mind, 
    \nomad is now able to manage advantageously a large number of cores.

The first release of \nomad4 will serve as a basis for future developments, starting with the integration of some of the improvements from the last 12 years:
\begin{itemize}
\item{\bimads~\cite{AuSaZg2008a}, \multimads~\cite{AuSaZg2010a} and
\dmultimads~\cite{BiLedSa2020} for multiobjective optimization;}
\item{\robustmads~\cite{AudIhaLedTrib2016} and \stomads~\cite{G-2019-30}
    for robust and stochastic optimization;}
\item{Variable Neighbourhood Search}~\cite{MlHa97a,AuBeLe08}
    to escape locally optimal solutions;
\item{Categorical}~\cite{AuDe01a} and {periodical variables}~\cite{AuLe2012}.
\end{itemize}

Another research direction is the application of the software to real industrial and engineering optimization problems
    to facilitate modeling, solving, analyzing and finding solutions for users.
Each application has its specificity, that may result in a new generic feature within \nomad.
Close collaborations with industry users is crucial for the development of \nomad. 
Contact us for projects, we look forward working with you.

\bigskip
{\bf Acknowledgments:}
This work is supported by the NSERC CRD RDCPJ~490744-15 grant coupled with an Innov\'E\'E grant, both in collaboration with Hydro-Qu\'ebec and Rio~Tinto, and by the NSERC Alliance grant 544900-19 in collaboration with Huawei-Canada.


\bibliographystyle{plain}
\bibliography{bibliography}

\begin{thebibliography}{10}

\bibitem{PhysRevD.87.042001}
J.~Aasi et~al.
\newblock {Einstein@Home all-sky search for periodic gravitational waves in
  LIGO S5 data}.
\newblock {\em Phys. Rev. D}, 87:042001, February 2013.

\bibitem{AACW09a}
M.A. Abramson, C.~Audet, J.W. Chrissis, and J.G. Walston.
\newblock {Mesh Adaptive Direct Search Algorithms for Mixed Variable
  Optimization}.
\newblock {\em Optimization Letters}, 3(1):35--47, 2009.

\bibitem{AbAuDeLe09}
M.A. Abramson, C.~Audet, J.E. {Dennis, Jr.}, and S.~{Le~Digabel}.
\newblock {OrthoMADS: A Deterministic MADS Instance with Orthogonal
  Directions}.
\newblock {\em SIAM Journal on Optimization}, 20(2):948--966, 2009.

\bibitem{AdAuBeYa2014}
L.~Adjengue, C.~Audet, and I.~Ben Yahia.
\newblock {A variance-based method to rank input variables of the Mesh Adaptive
  Direct Search algorithm}.
\newblock {\em Optimization Letters}, 8(5):1599--1610, 2014.

\bibitem{AlAuGhKoLed2020}
S.~Alarie, C.~Audet, A.E. Gheribi, M.~Kokkolaras, and S.~{Le~Digabel}.
\newblock {Two decades of blackbox optimization applications}.
\newblock Technical Report G-2020-58, Les cahiers du GERAD, 2020.

\bibitem{AmAuCoLed2016}
N.~Amaioua, C.~Audet, A.R. Conn, and S.~{Le~Digabel}.
\newblock {Efficient solution of quadratically constrained quadratic
  subproblems within a direct-search algorithm}.
\newblock {\em European Journal of Operational Research}, 268(1):13--24, 2018.

\bibitem{Audet2014a}
C.~Audet.
\newblock A survey on direct search methods for blackbox optimization and their
  applications.
\newblock In P.M. Pardalos and T.M. Rassias, editors, {\em Mathematics without
  boundaries: Surveys in interdisciplinary research}, chapter~2, pages 31--56.
  Springer, New York, NY, 2014.

\bibitem{AuBeCh2008a}
C.~Audet, V.~B\'echard, and J.~Chaouki.
\newblock Spent potliner treatment process optimization using a {MADS}
  algorithm.
\newblock {\em Optimization and Engineering}, 9(2):143--160, 2008.

\bibitem{AuBeLe08}
C.~Audet, V.~B\'echard, and S.~{Le~Digabel}.
\newblock {Nonsmooth optimization through Mesh Adaptive Direct Search and
  Variable Neighborhood Search}.
\newblock {\em Journal of Global Optimization}, 41(2):299--318, 2008.

\bibitem{AuCM2019}
C.~Audet and J.~C{\^o}t{\'e}-Massicotte.
\newblock {Dynamic improvements of static surrogates in direct search
  optimization}.
\newblock {\em Optimization Letters}, 13(6):1433--1447, 2019.

\bibitem{AuDe01a}
C.~Audet and J.E. {Dennis, Jr.}
\newblock Pattern search algorithms for mixed variable programming.
\newblock {\em SIAM Journal on Optimization}, 11(3):573--594, 2001.

\bibitem{AuDe2006}
C.~Audet and J.E. {Dennis, Jr.}
\newblock {Mesh Adaptive Direct Search Algorithms for Constrained
  Optimization}.
\newblock {\em SIAM Journal on Optimization}, 17(1):188--217, 2006.

\bibitem{AuDe09a}
C.~Audet and J.E. {Dennis, Jr.}
\newblock {A Progressive Barrier for Derivative-Free Nonlinear Programming}.
\newblock {\em SIAM Journal on Optimization}, 20(1):445--472, 2009.

\bibitem{AuDeLe07}
C.~Audet, J.E. {Dennis, Jr.}, and S.~{Le~Digabel}.
\newblock {Parallel Space Decomposition of the Mesh Adaptive Direct Search
  Algorithm}.
\newblock {\em SIAM Journal on Optimization}, 19(3):1150--1170, 2008.

\bibitem{AuLedDiSwMa2013}
C.~Audet, K.~Diest, S.~{Le~Digabel}, L.A. Sweatlock, and D.E. Marthaler.
\newblock {Metamaterial Design by Mesh Adaptive Direct Search}.
\newblock In {\em Numerical Methods for Metamaterial Design}, volume 127 of
  {\em Topics in Applied Physics}, pages 71--96. Springer, Dordrecht, 2013.

\bibitem{G-2019-30}
C.~Audet, K.J. Dzahini, M.~Kokkolaras, and S.~{Le~Digabel}.
\newblock {Stochastic mesh adaptive direct search for blackbox optimization
  using probabilistic estimates}.
\newblock Technical Report G-2019-30, Les cahiers du GERAD, 2021.
\newblock To appear in {\em Computational Optimization and Applications}.

\bibitem{AuHa2017}
C.~Audet and W.~Hare.
\newblock {\em {Derivative-Free and Blackbox Optimization}}.
\newblock Springer Series in Operations Research and Financial Engineering.
  Springer, Cham, Switzerland, 2017.

\bibitem{AuIaLeDTr2014}
C.~Audet, A.~Ianni, S.~{Le~Digabel}, and C.~Tribes.
\newblock {Reducing the Number of Function Evaluations in Mesh Adaptive Direct
  Search Algorithms}.
\newblock {\em SIAM Journal on Optimization}, 24(2):621--642, 2014.

\bibitem{AudIhaLedTrib2016}
C.~Audet, A.~Ihaddadene, S.~{Le~Digabel}, and C.~Tribes.
\newblock {Robust optimization of noisy blackbox problems using the Mesh
  Adaptive Direct Search algorithm}.
\newblock {\em Optimization Letters}, 12(4):675--689, 2018.

\bibitem{AuKoLedTa2016}
C.~Audet, M.~Kokkolaras, S.~{Le~Digabel}, and B.~Talgorn.
\newblock {Order-based error for managing ensembles of surrogates in mesh
  adaptive direct search}.
\newblock {\em Journal of Global Optimization}, 70(3):645--675, 2018.

\bibitem{AuLe2012}
C.~Audet and S.~{Le~Digabel}.
\newblock {The mesh adaptive direct search algorithm for periodic variables}.
\newblock {\em Pacific Journal of Optimization}, 8(1):103--119, 2012.

\bibitem{AuLedPey2014}
C.~Audet, S.~{Le~Digabel}, and M.~Peyrega.
\newblock {Linear equalities in blackbox optimization}.
\newblock {\em Computational Optimization and Applications}, 61(1):1--23, 2015.

\bibitem{AuLedTr2014}
C.~Audet, S.~{Le~Digabel}, and C.~Tribes.
\newblock {Dynamic scaling in the mesh adaptive direct search algorithm for
  blackbox optimization}.
\newblock {\em Optimization and Engineering}, 17(2):333--358, 2016.

\bibitem{AuLeDTr2018}
C.~Audet, S.~{Le~Digabel}, and C.~Tribes.
\newblock {The Mesh Adaptive Direct Search Algorithm for Granular and Discrete
  Variables}.
\newblock {\em SIAM Journal on Optimization}, 29(2):1164--1189, 2019.

\bibitem{AuSaZg2008a}
C.~Audet, G.~Savard, and W.~Zghal.
\newblock {Multiobjective Optimization Through a Series of Single-Objective
  Formulations}.
\newblock {\em SIAM Journal on Optimization}, 19(1):188--210, 2008.

\bibitem{AuSaZg2010a}
C.~Audet, G.~Savard, and W.~Zghal.
\newblock {A mesh adaptive direct search algorithm for multiobjective
  optimization}.
\newblock {\em European Journal of Operational Research}, 204(3):545--556,
  2010.

\bibitem{AuTr2018}
C.~Audet and C.~Tribes.
\newblock {Mesh-based Nelder-Mead algorithm for inequality constrained
  optimization}.
\newblock {\em Computational Optimization and Applications}, 71(2):331--352,
  2018.

\bibitem{BiLedSa2020}
J.~Bigeon, S.~{Le~Digabel}, and L.~Salomon.
\newblock {DMulti-MADS: Mesh adaptive direct multisearch for bound-constrained
  blackbox multiobjective optimization}.
\newblock Technical Report G-2020-25, Les cahiers du GERAD, 2021.
\newblock To appear in {\em Computational Optimization and Applications}.

\bibitem{ChWa2010}
X.~Chen and N.~Wang.
\newblock Optimization of short-time gasoline blending scheduling problem with
  a {DNA} based hybrid genetic algorithm.
\newblock {\em Chemical Engineering and Processing: Process Intensification},
  49(10):1076--1083, 2010.

\bibitem{CoLed2011}
A.R. Conn and S.~{Le~Digabel}.
\newblock Use of quadratic models with mesh-adaptive direct search for
  constrained black box optimization.
\newblock {\em Optimization Methods and Software}, 28(1):139--158, 2013.

\bibitem{FeMe1952}
E.~Fermi and N.~Metropolis.
\newblock Numerical solution of a minimum problem.
\newblock Los Alamos Unclassified Report LA--1492, Los Alamos National
  Laboratory, Los Alamos, USA, 1952.

\bibitem{MScMLG}
M.~Lemyre Garneau.
\newblock {Modelling of a solar thermal power plant for benchmarking blackbox
  optimization solvers}.
\newblock Master's thesis, Polytechnique Montr\'eal, 2015.
\newblock Code available at \url{https://github.com/bbopt/solar}.

\bibitem{GeHiLedAuTerScha2014}
E.M. Gertz, T.~Hiekkalinna, S.~{Le~Digabel}, C.~Audet, J.D. Terwilliger, and
  A.A. Schaffer.
\newblock {PSEUDOMARKER 2.0: efficient computation of likelihoods using NOMAD}.
\newblock {\em BMC Bioinformatics}, 15(47):8, 2014.

\bibitem{GhAuLeBeBaPe2012}
A.E. Gheribi, C.~Audet, S.~{Le~Digabel}, E.~B\'elisle, C.W. Bale, and A.D.
  Pelton.
\newblock {Calculating optimal conditions for alloy and process design using
  thermodynamic and properties databases, the FactSage software and the Mesh
  Adaptive Direct Search algorithm}.
\newblock {\em CALPHAD: Computer Coupling of Phase Diagrams and
  Thermochemistry}, 36:135--143, 2012.

\bibitem{GhHaBeRoChPeBaLe2014}
A.E. Gheribi, J.-P. Harvey, E.~B\'elisle, C.~Robelin, P.~Chartrand, A.D.
  Pelton, C.W. Bale, and S.~{Le~Digabel}.
\newblock {Use of a biobjective direct search algorithm in the process design
  of material science applications}.
\newblock {\em Optimization and Engineering}, 17(1):27--45, 2016.

\bibitem{GhLeAuCh2013}
A.E. Gheribi, S.~{Le~Digabel}, C.~Audet, and P.~Chartrand.
\newblock {Identifying optimal conditions for magnesium based alloy design
  using the Mesh Adaptive Direct Search algorithm}.
\newblock {\em Thermochimica Acta}, 559(0):107--110, 2013.

\bibitem{GhRoLeAuPe2011}
A.E. Gheribi, C.~Robelin, S.~{Le~Digabel}, C.~Audet, and A.D. Pelton.
\newblock {Calculating all local minima on liquidus surfaces using the FactSage
  software and databases and the Mesh Adaptive Direct Search algorithm}.
\newblock {\em The Journal of Chemical Thermodynamics}, 43(9):1323--1330, 2011.

\bibitem{GoOrTo03}
{N.I.M.} Gould, D.~Orban, and {Ph.}L. Toint.
\newblock {CUTEr (and SifDec): A Constrained and Unconstrained Testing
  Environment, revisited}.
\newblock {\em ACM Transactions on Mathematical Software}, 29(4):373--394,
  2003.

\bibitem{HaMl01a}
P.~Hansen and N.~Mladenovi\'c.
\newblock Variable neighborhood search: principles and applications.
\newblock {\em European Journal of Operational Research}, 130(3):449--467,
  2001.

\bibitem{HaBeAuKo03a}
R.E. Hayes, F.H. Bertrand, C.~Audet, and S.T. Kolaczkowski.
\newblock Catalytic combustion kinetics: Using a direct search algorithm to
  evaluate kinetic parameters from light-off curves.
\newblock {\em The Canadian Journal of Chemical Engineering}, 81(6):1192--1199,
  2003.

\bibitem{HiHeEu2020}
A.~Hibberd, A.M. Hein, and T.M. Eubanks.
\newblock {Project Lyra: Catching 1I/{\textquoteleft}Oumuamua - Mission
  opportunities after 2024}.
\newblock {\em Acta Astronautica}, 170:136--144, 2020.

\bibitem{HoSc1981}
W.~Hock and K.~Schittkowski.
\newblock {\em {Test Examples for Nonlinear Programming Codes}}, volume 187 of
  {\em Lecture Notes in Economics and Mathematical Systems}.
\newblock Springer, Berlin, Germany, 1981.

\bibitem{KiArYa2011}
S.~Kitayama, M.~Arakawa, and K.~Yamazaki.
\newblock Sequential approximate optimization using radial basis function
  network for engineering optimization.
\newblock {\em Optimization and Engineering}, 12(4):535--557, 2011.

\bibitem{MScGL}
G.~Lameynardie.
\newblock {Sondes locales intensives lors de l'ex\'ecution de l'algorithme MADS
  dans un environnement parall\`ele}.
\newblock Master's thesis, Polytechnique Montr\'eal, 2020.

\bibitem{Le09b}
S.~{Le~Digabel}.
\newblock {Algorithm 909: NOMAD: Nonlinear Optimization with the MADS
  algorithm}.
\newblock {\em {ACM} Transactions on Mathematical Software}, 37(4):44:1--44:15,
  2011.

\bibitem{LiTr2017}
G.~Liuzzi and K.~Truemper.
\newblock {Parallelized hybrid optimization methods for nonsmooth problems
  using NOMAD and linesearch}.
\newblock {\em Computational and Applied Mathematics}, 37(3):3172--3207, 2018.

\bibitem{LuVl00}
L.~Luk{\v s}an and J.~Vl{\v c}ek.
\newblock Test problems for nonsmooth unconstrained and linearly constrained
  optimization.
\newblock Technical Report V-798, ICS AS CR, 2000.

\bibitem{McCoBe79a}
M.D. McKay, R.J. Beckman, and W.J. Conover.
\newblock A comparison of three methods for selecting values of input variables
  in the analysis of output from a computer code.
\newblock {\em Technometrics}, 21(2):239--245, 1979.

\bibitem{MTSMEBKF2019}
K.~Mehrgan, J.~Thomas, R.~Saglia, X.~Mazzalay, P.~Erwin, R.~Bender, M.~Kluge,
  and M.~Fabricius.
\newblock {A 40 Billion Solar-mass Black Hole in the Extreme Core of Holm 15A,
  the Central Galaxy of Abell 85}.
\newblock {\em The Astrophysical Journal}, 887(2):195, 2019.

\bibitem{MezCoe05}
E.~Mezura-Montes and C.A. Coello.
\newblock {Useful Infeasible Solutions in Engineering Optimization with
  Evolutionary Algorithms}.
\newblock In {\em {Proceedings of the 4th Mexican International Conference on
  Advances in Artificial Intelligence}}, MICAI'05, pages 652--662, Berlin,
  Heidelberg, 2005. Springer-Verlag.

\bibitem{MiCaGuLeAuLe2014}
M.~Minville, D.~Cartier, C.~Guay, L.-A. Leclaire, C.~Audet, S.~{Le~Digabel},
  and J.~Merleau.
\newblock {Improving process representation in conceptual hydrological model
  calibration using climate simulations}.
\newblock {\em Water Resources Research}, 50:5044--5073, 2014.

\bibitem{MlHa97a}
N.~Mladenovi\'c and P.~Hansen.
\newblock Variable neighborhood search.
\newblock {\em Computers and Operations Research}, 24(11):1097--1100, 1997.

\bibitem{MoWi2009}
J.J. Mor\'e and S.M. Wild.
\newblock {Benchmarking Derivative-Free Optimization Algorithms}.
\newblock {\em SIAM Journal on Optimization}, 20(1):172--191, 2009.

\bibitem{RNFSSSRGMR2019}
V.~Ramakrishnan, N.M. Nagar, C.~Finlez, T.~Storchi-Bergmann, R.~Slater,
  A.~Schnorr-Muller, R.A. Riffel, C.G. Mundell, and A.~Robinson.
\newblock {Nuclear kinematics in nearby AGN - I. An ALMA perspective on the
  morphology and kinematics of the molecular CO(2–1) emission}.
\newblock {\em Monthly Notices of the Royal Astronomical Society},
  487(1):444--455, 2019.

\bibitem{RodRenWat98}
J.F. Rodr{\'\i}guez, J.E. Renaud, and L.T. Watson.
\newblock {Trust Region Augmented Lagrangian Methods for Sequential Response
  Surface Approximation and Optimization}.
\newblock {\em Journal of Mechanical Design}, 120(1):58--66, 1998.

\bibitem{MScLASMC}
L.A. {Sarrazin-Mc~Cann}.
\newblock {Optimisation et ordonnancement en optimisation sans d\'eriv\'ees}.
\newblock Master's thesis, Polytechnique Montr\'eal, 2018.

\bibitem{SeCoAu2017}
S.~S\'eguin, P.~C\^ot\'e, and C.~Audet.
\newblock {Scenario tree modeling for stochastic short-term hydropower
  operations planning}.
\newblock {\em Journal of Water Resources Planning and Management},
  143(12):04017073--1--12, 2017.

\bibitem{sgtelib}
B.~Talgorn.
\newblock {{\sf sgtelib}: Surrogate model library for Derivative-Free
  Optimization}, 2019.

\bibitem{TaAuKoLed2016}
B.~Talgorn, C.~Audet, M.~Kokkolaras, and S.~{Le~Digabel}.
\newblock {Locally weighted regression models for surrogate-assisted design
  optimization}.
\newblock {\em Optimization and Engineering}, 19(1):213--238, 2018.

\bibitem{TaLeDKo2014}
B.~Talgorn, S.~{Le~Digabel}, and M.~Kokkolaras.
\newblock {Statistical Surrogate Formulations for Simulation-Based Design
  Optimization}.
\newblock {\em Journal of Mechanical Design}, 137(2):021405--1--021405--18,
  2015.

\bibitem{Tang93a}
B.~Tang.
\newblock Orthogonal array-based latin hypercubes.
\newblock {\em Journal of the American Statistical Association},
  88(424):1392--1397, 1993.

\bibitem{TaoWan08}
J.~Tao and N.~Wang.
\newblock {DNA Double Helix Based Hybrid GA for the Gasoline Blending Recipe
  Optimization Problem}.
\newblock {\em Chemical Engineering and Technology}, 31(3):440--451, 2008.

\bibitem{Torc97a}
V.~Torczon.
\newblock On the convergence of pattern search algorithms.
\newblock {\em SIAM Journal on Optimization}, 7(1):1--25, 1997.

\bibitem{TriDuTre04a}
C.~Tribes, J.-F. Dub\'e, and J.-Y. Tr\'epanier.
\newblock Decomposition of multidisciplinary optimization problems:
  formulations and application to a simplified wing design.
\newblock {\em Engineering Optimization}, 37(8):775--796, 2005.

\bibitem{ZhaWan2010b}
J.~Zhao and N.~Wang.
\newblock {A bio-inspired algorithm based on membrane computing and its
  application to gasoline blending scheduling}.
\newblock {\em Computers and Chemical Engineering}, 35(2):272--283, 2011.

\end{thebibliography}

\end{document}